\documentclass[12pt,a4paper]{article}
\usepackage{amsfonts}
\usepackage{amsmath}
\usepackage{amssymb}
\usepackage[english]{babel}
\usepackage[T1]{fontenc}
\usepackage{graphicx}
\usepackage{subfigure}
\newtheorem{prop}{Proposition}[section]
\newtheorem{rem}[prop]{Remark}

\newtheorem{theo}[prop]{Theorem}
\newtheorem{cor}[prop]{Corollary}

\usepackage{epsfig}
\usepackage{latexsym}
\usepackage{amstext}
\usepackage{epsfig, graphicx}
\usepackage{epsfig}
\usepackage{latexsym}
\usepackage{amstext}
\usepackage{amssymb}
\usepackage{graphics}
\usepackage{enumerate}
\usepackage[]{natbib}

\usepackage{fullpage}

\newcommand{\beq}{\begin{eqnarray}}
\newcommand{\beqq}{\begin{eqnarray*}}
\newcommand{\eeq}{\end{eqnarray}}
\newcommand{\eeqq}{\end{eqnarray*}}

\def\QED{\quad\hbox{\hskip 4pt\vrule width 5pt height 6pt depth 1.5pt}}

\title{A Central Limit Theorem for a sequence of Brownian motions in the unit sphere in $\mathbb{R}^{n}$}
\author{ S. Vakeroudis\thanks{Laboratoire de Probabilit\'{e}s et Mod\`{e}les
Al\'{e}atoires (LPMA) CNRS : UMR7599,  Universit\'{e} Pierre et Marie
Curie - Paris VI,  Universit\'{e} Paris-Diderot - Paris VII, 4 place Jussieu, 75252 Paris Cedex 05, France.
E-mail: stavros.vakeroudis@etu.upmc.fr }
\thanks{Probability and Statistics Group, School of Mathematics, University of Manchester,
Alan Turing Building, Oxford Road, Manchester M13 9PL, United Kingdom. } \
M. Yor$^{\ast}$\thanks{Institut Universitaire de
France, Paris, France.} }
\date{\today}
\begin{document}

\maketitle
\begin{abstract}
We use a Stochastic Differential Equation satisfied by Brownian motion taking values in the unit sphere $S_{n-1}\subset\mathbb{R}^{n}$
and we obtain a Central Limit Theorem for a sequence of such Brownian motions.
We also generalize the results to the case of the $n$-dimensional Ornstein-Uhlenbeck processes.
\end{abstract}


$\vspace{5pt}$

\section{Introduction}
\renewcommand{\thefootnote}{\fnsymbol{footnote}}
This paper may be regarded as an extension to higher dimensions of the 2-dimensional study made in \cite{VYH11}. \\
We now consider a sequence of Brownian motions $\left(\Theta^{(k)}_{t},t\geq0\right)$, $k\in\mathbb{N}$ taking values in the unit sphere $S_{n-1}(\subset\mathbb{R}^{n})$,
all starting from the same point on the sphere. In Section \ref{secBM}, we introduce a general
representation of $\Theta^{(k)}$ in terms of a Stochastic Differential Equation. Using this representation, we describe
in detail in Section \ref{secTCL} the limit in law, as $K\rightarrow\infty$, for the renormalized sum:
\beqq
Z^{K}_{t}\equiv\frac{1}{\sqrt{K}} \sum^{K}_{k=1} \left(\Theta^{(k)}_{t}-E\left[\Theta^{(k)}_{t}\right]\right)
\eeqq
of these processes, indexed by $t\geq0$, and taking values in $\mathbb{R}^{n}$. Of course, one could
invoke the classical Central Limit Theorem (CLT), at least for the finite dimensional marginals of
$\left(Z^{K}_{t},t\geq0\right)$, as $K\rightarrow\infty$. However, with the help of stochastic calculus,
there is much more to say about the description of the asymptotics.
Finally, in this Section, we remark that the CLT can be generalized to Ornstein-Uhlenbeck processes
taking values in the unit sphere $S_{n-1}(\subset\mathbb{R}^{n})$.
Three technical points are gathered in an Appendix.

Further extensions may also be obtained, by following e.g. \cite{Ito83} or \cite{Ochi85} and studying for which class of functions $f(\Theta)$ we can obtain a functional CLT such as (\ref{Masymp}) (see below) for $f(\Theta_{t})$, instead of the unique function
$f_{*}(\Theta)=\Theta$ which we study here.

\section{A presentation of Brownian motion in the sphere $S_{n-1}$} \label{secBM}

As remarked in \cite{Str71} and \cite{Yor84} (eq. (4.j), p.34),
Brownian motion $\left(\Theta_{t},t\geq0\right)$
in the unit sphere $S_{n-1}\subset\mathbb{R}^{n}$ may be viewed as the solution of a Stochastic Differential Equation:
\beq\label{theta}
\Theta_{t}=\Theta_{0}+\int^{t}_{0}\sigma^{0,1}(\Theta_{s})\cdot dB_{s}-\frac{n-1}{2}\int^{t}_{0}ds \ \Theta_{s} \ .
\eeq
In (\ref{theta}), $B_{t}\equiv\left(B^{(i)}_{t},i\leq n\right)$, $t\geq0$, denotes a $n$-dimensional Brownian
motion starting from $a\neq0$,
while $\left(\sigma^{0,1}(x),x\in S_{n-1}\right)$ denotes the family of $n\times n$ matrices (see e.g. \cite{Kry80}), defined by:
\beq\label{sigma}
\sigma^{0,1}(x)=\left(\delta_{i,j}-x_{i}x_{j}\right)_{i,j\leq n}, \ \ \ (x\in S_{n-1})
\eeq
and/or characterized by:
\beq\label{sigma2}
\sigma^{0,1}(x)\cdot x=0, \ \ \textrm{and} \ \ \sigma^{0,1}(x)\cdot y=y, \ \ \textrm{if} \ y\cdot x=0.
\eeq
Note that $\sigma^{0,1}(x)$ is symmetric and satisfies: $\sigma^{0,1}(x)\sigma^{0,1}(x)=\sigma^{0,1}(x)$. \\
Thus, from (\ref{sigma2}), we deduce that:
\beq
\sigma^{0,1}(x)m &=& m-\left(m\cdot x\right)x, \ m\in\mathbb{R}^{n}, \label{sigmasym} \\
\left(\sigma^{0,1}(x)m\right)\cdot\left(\sigma^{0,1}(x)m'\right)&=&\left(m\cdot m'\right)-\left(m\cdot x\right)\left(m'\cdot x\right), m,m'\in\mathbb{R}^{n}. \label{sigmasymbis}
\eeq

\section{A Central Limit Theorem for a sequence of Brownian motions in the sphere $S_{n-1}(\subset\mathbb{R}^{n})$} \label{secTCL}

Let $\Theta^{(1)},\ldots,\Theta^{(k)},\ldots$ be a sequence of such independent and identically distributed Brownian motions
in the sphere $S_{n-1}$. We aim for a Central Limit Theorem concerning:
\beq
Z^{K}_{t}\equiv\frac{1}{\sqrt{K}} \sum^{K}_{k=1} \left(\Theta^{(k)}_{t}-E\left[\Theta^{(k)}_{t}\right]\right).
\eeq
Adding $K$ equations of the kind of (\ref{theta}) term by term, for $\left(\Theta^{(k)}_{t},k\leq K\right)$, it is immediate that:
\beq\label{Z}
Z^{K}_{t}=M^{K}_{t}-\frac{n-1}{2} \int^{t}_{0} ds \; Z^{K}_{s},
\eeq
with
\beq
M^{K}_{t}=\frac{1}{\sqrt{K}} \sum^{K}_{k=1} \int^{t}_{0} \sigma^{0,1}(\Theta^{(k)}_{s}) \cdot dB^{(k)}_{s}.
\eeq
Thus, from (\ref{Z}), we obtain:
\beq\label{ZMK}
Z^{K}_{t}=\exp\left(\frac{-(n-1)t}{2}\right) \int^{t}_{0} \exp\left(\frac{(n-1)s}{2}\right) \; dM^{K}_{s}.
\eeq
Now, clearly, the Central Limit Theorem for $\left(Z^{K}_{t}\right)$, $K\rightarrow\infty$, which we are seeking,
will follow from the limit in law of the martingales $\left(M^{K}_{t},t\geq0\right)$, as $K\rightarrow\infty$.
We now state both limit results in the following:
\begin{theo}\label{theoasymp}
$\left.a\right)$ The sequence of martingales $\left(M^{K}_{t},t\geq0\right)$ converges in law, as $K\rightarrow\infty$,
towards:
\beq\label{Minf3}
M^{(\infty)}_{t} &=& \sqrt{1-\frac{1}{n}}\left\{\Theta(0) \int^{t}_{0} \sqrt{1-e^{-ns}} d\beta_{s} +\int^{t}_{0} \sqrt{1+\frac{e^{-ns}}{n-1}}dB'_{s}  \right\},
\eeq
where $\left(\beta_{s},s\geq0\right)$ is a 1-dimensional BM and $\left(B'_{s},s\geq0\right)$ is a $(n-1)-$dimensional BM
taking values in the hyperplane which is orthogonal to $\Theta(0)$, and $B'$ is independent of $\beta$. \\
$\left.b\right)$ Consequently, $\left(Z^{K}_{t},t\geq0\right)$ converges in law, as $K\rightarrow\infty$,
towards:
\beq\label{ZMinf}
Z^{(\infty)}_{t}=\exp\left(-\frac{n-1}{2}t\right) \int^{t}_{0}\exp\left(\frac{n-1}{2}s\right) dM^{(\infty)}_{s}.
\eeq
\end{theo}
{\noindent \textbf{Proof of Theorem \ref{theoasymp}:} }
Using the Law of Large Numbers, it is not difficult to show that\footnote[8]{In the Appendix \ref{apgen}, a more general result, concerning $\frac{1}{\sqrt{K}}\int^{t}_{0}\sum^{K}_{k=1}H^{(k)}_{s}\cdot dB^{(k)}_{s}$ is presented, where $(B^{(k)}, k=1,\ldots,K)$ are $K$ independent BMs
and $(H^{(k)},B^{(k)})$ are $k\leq K$ iid random vectors.}:
\beq\label{Masymp}
\left(M^{K}_{t},t\geq0\right)&\overset{{(law)}}{\underset{K\rightarrow\infty}\longrightarrow}&
\int^{t}_{0}\sqrt{E\left[\sigma^{0,1}(\Theta^{(1)}_{s})\right]}\cdot dB^{(1)}_{s}\equiv\left(M^{\infty}_{t},t\geq0\right),
\eeq
where $Q(s)\equiv E\left[\sigma^{0,1}(\Theta^{(1)}_{s})\right]$ is a deterministic matrix,
depending on $s$. The RHS of (\ref{Masymp}) is a centered Gaussian martingale in $\mathbb{R}^{n}$.
Before computing the square root involved in (\ref{Masymp}), we shall first calculate (see (\ref{sigma}) for the definition of $\sigma^{0,1}$):
\beq
E\left[\sigma^{0,1}(\Theta^{(1)}_{s})\right]=\left(\delta_{i,j}-E\left[\Theta_{i}(s)\Theta_{j}(s)\right]\right)_{i,j\leq n}.
\eeq
In order to calculate $E\left[\Theta^{(i)}_{s}\Theta^{(j)}_{s}\right]$ as "naturally" as possible,
we consider two generic vectors $m$ and $m'$ in $\mathbb{R}^{n}$, and we compute:
\beq
\varphi_{m,m'}(t)\equiv E\left[\left(m \cdot \Theta_{t}\right)\left(m' \cdot\Theta_{t}\right)\right].
\eeq
Using (\ref{theta}) and the (special) properties of the matrices $\left\{\sigma^{0,1}(x)\right\}$, we easily
deduce from It\^{o}'s formula, that:
\beq\label{mtheta}
E\left[\left(m \cdot \Theta_{t}\right)\left(m' \cdot\Theta_{t}\right)\right]&=&\left(m \cdot \Theta_{0}\right)\left(m' \cdot\Theta_{0}\right)
- (n-1) \int^{t}_{0}ds \; E\left[\left(m \cdot \Theta_{s}\right)\left(m' \cdot\Theta_{s}\right)\right] \nonumber \\
&& \ \ \ \ + \int^{t}_{0}ds \; E\left[\left(\sigma^{0,1}(\Theta_{s})m \right)\cdot\left(\sigma^{0,1}(\Theta_{s})m'\right)\right] .
\eeq
Using (\ref{sigmasymbis}), (\ref{mtheta}) simplifies as:
\beq\label{mtheta2}
&& E\left[\left(m \cdot \Theta_{t}\right)\left(m' \cdot\Theta_{t}\right)\right]=\left(m \cdot \Theta_{0}\right)\left(m' \cdot\Theta_{0}\right)
- (n-1) \int^{t}_{0}ds \; E\left[\left(m \cdot \Theta_{s}\right)\left(m' \cdot\Theta_{s}\right)\right] \nonumber \\
&& \ \ \ \ \ \ \ \ \ \ \ \ \ \ \ \ \ \ \ \ \ \ \ \ \ \ \ \ \ \ \ \ \ \ \ \ + \int^{t}_{0}ds \; \left(m\cdot m'-E\left[\left(m \cdot \Theta_{s}\right)\left(m' \cdot\Theta_{s}\right)\right]\right) \nonumber \\
&& = \ \ \ \ \left(m \cdot \Theta_{0}\right)\left(m' \cdot\Theta_{0}\right)+\left(m\cdot m'\right)t-n\int^{t}_{0}ds \; E\left[\left(m \cdot \Theta_{s}\right)\left(m' \cdot\Theta_{s}\right)\right].
\eeq
Consequently, the function $\varphi_{m,m'}(t)=E\left[\left(m \cdot \Theta_{t}\right)\left(m' \cdot\Theta_{t}\right)\right]$
is the solution of a first order linear differential equation, hence:
\beq
E\left[\left(m \cdot \Theta_{t}\right)\left(m' \cdot\Theta_{t}\right)\right]=e^{-nt} \left\{\left(m \cdot \Theta_{0}\right)\left(m' \cdot\Theta_{0}\right)+\left(m\cdot m'\right)\int^{t}_{0}e^{ns}ds \right\}.
\eeq
Now, taking $m=e_{i}$ and $m'=e_{j}$, where $\left(e_{k};k\leq n\right)$ is the canonical basis of $\mathbb{R}^{n}$,
the matrix $Q(s)$ has elements:
\beq
\mathrm{for} \ i\neq j, \ \ \left(Q(s)\right)_{i,j}=-E\left[\Theta_{i}(s)\Theta_{j}(s)\right]=-\Theta_{i}(0)\Theta_{j}(0)e^{-ns},
\eeq
\beq
\mathrm{for} \ i=j, \ \ \left(Q(s)\right)_{i,i}&=&1-E\left[\Theta_{i}(s)\Theta_{i}(s)\right] \nonumber \\
&=&1-\left\{\left(\Theta_{i}(0)\right)^{2}e^{-ns}+e^{-ns}\left(\frac{e^{ns}-1}{n}\right)\right\} \nonumber \\
&=&\left(1-\frac{1}{n}\right)+e^{-ns}\left(\frac{1}{n}-\left(\Theta_{i}(0)\right)^{2}\right).
\eeq
Finally:
\beq\label{Q}
Q(s)&=&\left(1-\frac{1}{n}\right)Id+e^{-ns}\left(\frac{1}{n}\delta_{ij}-\Theta_{i}(0)\Theta_{j}(0)\right)_{i,j\leq n} \nonumber \\
&\equiv& \left(1-\frac{1}{n}\left(1-e^{-ns}\right)\right)Id-e^{-ns}\left(\Theta_{i}(0)\Theta_{j}(0)\right)_{i,j\leq n} \ .
\eeq
Using (\ref{Q}) in the RHS of (\ref{Masymp}), we obtain:
\beq\label{Minf}
M^{(\infty)}_{t}\equiv\int^{t}_{0} \sqrt{Q(s)} \; dB_{s}.
\eeq
Now $\sqrt{Q(s)}\equiv\Lambda(s)$, where (for the explicit calculation, see Appendix \ref{aproot}):
\beq
\Lambda(s)&\equiv&\sqrt{1-\frac{1}{n}}\sqrt{1-e^{-ns}} Id+\sqrt{1-\frac{1}{n}}\left(\sqrt{1+\frac{e^{-ns}}{n-1}}-\sqrt{1-e^{-ns}}\right)
\sigma^{0,1}\left(\Theta(0)\right) \nonumber \\
&=& \sqrt{1-\frac{1}{n}}\left\{\sqrt{1-e^{-ns}} Id+
\left(\sqrt{1+\frac{e^{-ns}}{n-1}}-\sqrt{1-e^{-ns}}\right)\sigma^{0,1}\left(\Theta(0)\right)\right\}.
\eeq
Thus, (\ref{Minf}) now writes:
\beq\label{Minf2}
M^{(\infty)}_{t}&\equiv&\int^{t}_{0} \Lambda(s) \; dB_{s} \nonumber \\
&=& \sqrt{1-\frac{1}{n}}\left\{\int^{t}_{0} \sqrt{1-e^{-ns}} dB_{s} +\int^{t}_{0} \left[\sqrt{1+\frac{e^{-ns}}{n-1}}-\sqrt{1-e^{-ns}}\right]
\sigma^{0,1}\left(\Theta(0)\right)dB_{s} \right\}. \nonumber \\
\eeq
We remark here that, with $\beta_{s}\equiv \Theta(0)\cdot B_{s}$,
\beq\label{Bprime}
B'_{s}=B_{s}-\Theta(0)\beta_{s}\equiv\sigma^{0,1}\left(\Theta(0)\right)B_{s}
\eeq
is a $(n-1)$-dimensional BM taking values in the hyperplane which is orthogonal to $\Theta(0)$.
Thus, from (\ref{Minf2}) we deduce (\ref{Minf3}). \\
From (\ref{ZMK}), letting $K\rightarrow\infty$, we obtain (\ref{ZMinf}).
\hfill \QED
\\ \\
Moreover, changing the variables $s=t-u$ and using the dominated convergence Theorem, we have:
\beq\label{ZMinf2}
Z^{(\infty)}_{t} &\stackrel{(law)}{=}& \sqrt{1-\frac{1}{n}} \exp\left(-\frac{n-1}{2}t\right) \int^{t}_{0}
\exp\left(\frac{n-1}{2}s\right) \times \nonumber \\
&& \ \ \ \ \ \ \ \ \ \ \ \ \ \ \ \ \ \ \ \ \ \ \ \ \ \ \ \ \ \times \left\{\Theta(0) \sqrt{1-e^{-ns}} d\beta_{s}+\sqrt{1+\frac{e^{-ns}}{n-1}}dB'_{s}  \right\} \label{ZMinf2bis} \\
&\overset{{s=t-u}}{\underset{(law)}=}& \sqrt{1-\frac{1}{n}} \int^{t}_{0} \exp\left(-\frac{n-1}{2}u\right) \left\{\sqrt{1-e^{-n(t-u)}} \Theta(0) d\beta_{u}+\sqrt{1+\frac{e^{-n(t-u)}}{n-1}} dB'_{u}\right\} \nonumber \\
&\overset{{t\rightarrow\infty}}{\longrightarrow}& \sqrt{1-\frac{1}{n}} \int^{\infty}_{0} \exp\left(-\frac{n-1}{2}u\right) \Theta(0) d\beta_{u} +\sqrt{1-\frac{1}{n}} \int^{\infty}_{0} \exp\left(-\frac{n-1}{2}u\right) dB'_{u}. \nonumber \\
\eeq
\begin{prop}\label{propZMinfL}
The following asymptotic results hold: \\
$\left.a\right)$
\beq
Z^{(\infty)}_{t} \overset{{(law)}}{\underset{t\rightarrow\infty}\longrightarrow} Z^{(\infty)}_{\infty},
\eeq
where:
\beq\label{Zinfinf}
Z^{(\infty)}_{\infty}\equiv \sqrt{1-\frac{1}{n}} \int^{\infty}_{0} \exp\left(-\frac{n-1}{2}u\right) dB_{u} \ .
\eeq
$\left.b\right)$
\beq\label{ZMinf3}
Z^{(\infty)}_{t}-\exp\left(-\frac{n-1}{2}t\right) \int^{t}_{0} \sqrt{1-\frac{1}{n}} \exp\left(\frac{n-1}{2}s\right) dB_{s} \overset{{L^{2}}}{\underset{t\rightarrow\infty}\longrightarrow} 0.
\eeq
\end{prop}
Part $\left.a\right)$ of Proposition \ref{propZMinfL} follows from the previous calculations, using (\ref{Bprime}).
In order to prove part $\left.b\right)$, it suffices to use the expression (\ref{ZMinf2bis})
and the following Proposition, which reinforces the convergence in $L^{2}$
result in (\ref{ZMinf3}).
\begin{prop}\label{propositionZ}
As $t\rightarrow\infty$, the Gaussian martingales:
\beq
    \left(G^{(0)}_{t}\Theta(0),t\geq0\right)\equiv\Theta(0)\left(\int^{t}_{0} \sqrt{1-e^{-ns}} e^{\frac{n-1}{2}s} d\beta_{s}
-\int^{t}_{0} e^{\frac{n-1}{2}s} d\beta_{s}, t\geq0\right),
\eeq
and
\beq
    \left(G'_{t},t\geq0\right)\equiv\left(\int^{t}_{0} \sqrt{1+\frac{e^{-ns}}{n-1}}e^{\frac{n-1}{2}s}dB'_{s}-\int^{t}_{0} e^{\frac{n-1}{2}s} \; dB'_{s}, t\geq0\right)
\eeq
converge a.s. and in $L^{2}$, and the limit variables are Gaussian, with variances, respectively:
$\left(\frac{\sqrt{\pi}\Gamma(-1+\frac{1}{n})}{n\Gamma(\frac{1}{2}+\frac{1}{n})}-\frac{n+1}{n-1}\right)$,
and $\frac{2 \; _{2}F_{1}\left(-\frac{1}{2},-1+\frac{1}{n},\frac{1}{n},\frac{1}{1-n}\right)-1}{n-1}$.
\end{prop}
{\noindent \textbf{Proof of Proposition \ref{propositionZ}}:} \\
$\left.a\right)$ The increasing process of the real-valued Gaussian martingale $G^{(0)}_{t}$ is:
\beqq
    \int^{t}_{0} e^{(n-1)s}\left(\sqrt{1-e^{-ns}}-1\right)^{2} ds,
\eeqq
which converges, as $t\rightarrow\infty$; thus:
\beqq
G^{(0)}_{t}{\underset{t\rightarrow\infty}\longrightarrow}\int^{\infty}_{0} \left(\sqrt{1-e^{-ns}} e^{\frac{n-1}{2}s}-e^{\frac{n-1}{2}s}\right) d\beta_{s},
\eeqq
where the convergence holds both a.s. and in every $L^{p}$. Of course, the limit variable is
Gaussian and its variance is given by (we change the variables $u=e^{-ns}$ and $B(a,b)$ denotes the Beta function with arguments $a$ and $b$\footnote[7]{We recall that $\left(\Gamma(x),x\geq0\right)$ denotes the Gamma function and $B(a,b)=\frac{\Gamma(a)\Gamma(b)}{\Gamma(a+b)}$.}):
\beqq
    && \int^{\infty}_{0} ds \; e^{(n-1)s} \left(\sqrt{1-e^{-ns}} - 1\right)^{2} = \frac{1}{n}\int^{1}_{0} du \; u^{-2+\frac{1}{n}} \left(\sqrt{1-u} - 1\right)^{2} \\
    &=& \frac{1}{n} \left[ \int^{1}_{0} du \; u^{-2+\frac{1}{n}} \left((1-u)-2\sqrt{1-u}+1\right)\right] \\
    &=& \frac{1}{n} \left\{B\left(-1+\frac{1}{n},2\right)-2B\left(-1+\frac{1}{n},\frac{3}{2}\right)-\frac{n}{n-1}\right\} \\
    &=& \frac{\sqrt{\pi}\Gamma(-1+\frac{1}{n})}{n\Gamma(\frac{1}{2}+\frac{1}{n})}-\frac{n+1}{n-1} \ .
\eeqq
To be rigorous, the integral $\int^{1}_{0} du \; u^{-\alpha} \left(\sqrt{1-u} - 1\right)^{2}$, which is well defined for $0<\alpha<1$,
can be extended analytically for any complex $\alpha$ with $\mathrm{Re}(\alpha)<3$. \\
$\left.b\right)$ Likewise, the "increasing process" of the vector-valued Gaussian martingale $G'_{t}$ is:
\beqq
    \int^{t}_{0} e^{(n-1)s}\left(\sqrt{1+\frac{e^{-ns}}{n-1}}-1\right)^{2} ds
\eeqq
which also converges as $t\rightarrow\infty$. The limit variable:
\beqq
\int^{\infty}_{0} \left(\sqrt{1+\frac{e^{-ns}}{n-1}}e^{\frac{n-1}{2}s}- e^{\frac{n-1}{2}s}\right) dB'_{s},
\eeqq
is also Gaussian and, by repeating the previous calculation, we easily compute its variance.
\begin{flushright}
\QED
\end{flushright}

{\noindent{\textbf{Proof of Proposition \ref{propZMinfL}: }}} \\
From Proposition \ref{propositionZ}, by multiplying both processes by $e^{\left(-\frac{(n-1)t}{2}\right)}\sqrt{1-\frac{1}{n}}$,
we obtain (\ref{ZMinf3}).
\begin{flushright}
\QED
\end{flushright}

\begin{rem} \textbf{(The Ornstein-Uhlenbeck case)} \\
In fact, for every process satisfying:
\beq
    dZ_{s} = d\mathbb{B}_{s} + h(|Z_{s}|) Z_{s} ds,
\eeq
where $\left(\mathbb{B}_{t},t\geq0\right)$ is a $n$-dimensional Brownian motion
(BM) and $h:\mathbb{R}_{+}\rightarrow\mathbb{R}$ is a bounded function, there is a CLT of the kind of Theorem \ref{theoasymp}. See Appendix \ref{secOU} for the proof.
\end{rem}

\appendix
\section{Appendix}

\subsection{Generalization for a class of symmetric matrices}\label{apgen}
For $K$ independent Brownian motions, and a class of symmetric matrices $H^{(k)}$ such that $(H^{(k)},B^{(k)})_{k\leq K}$ are iid,
we have:
\beq
\tilde{M}^{(K)}_{t}&\equiv& \frac{1}{\sqrt{K}}\int^{t}_{0}\sum^{K}_{k=1}H^{(k)}_{s}\cdot dB^{(k)}_{s}\overset{{(law)}}{\underset{K\rightarrow\infty}\longrightarrow} \int^{t}_{0} h_{s}dB_{s},
\eeq
with $h_{s}$ a deterministic symmetric positive definite matrix and $(B_{t},t\geq0)$ a $n$-dimensional BM. \\
Indeed, using $m$ a generic vector in $\mathbb{R}^{n}$, we have:
\beq
m\cdot\tilde{M}^{(K)}_{t}&\equiv& m\cdot\frac{1}{\sqrt{K}}\int^{t}_{0}\sum^{K}_{k=1}H^{(k)}_{s}\cdot dB^{(k)}_{s} = \frac{1}{\sqrt{K}}\sum^{K}_{k=1}\int^{t}_{0}\left(H^{(k)}_{s}m\right)\cdot dB^{(k)}_{s} \nonumber \\
&\overset{{(law)}}{\underset{K\rightarrow\infty}\longrightarrow}& \mathcal{N} \left(0; \int^{t}_{0} ds E\left[|H^{1}_{s}\cdot m|^{2}\right]\right),
\eeq
and for the variance, we have:
\beq
\int^{t}_{0} ds E\left[|H^{1}_{s}\cdot m|^{2}\right]=m\cdot \int^{t}_{0} ds E\left[H^{1}_{s}\tilde{H}^{1}_{s}\right]m,
\eeq
and
\beq
E\left[H^{1}_{s}\right]\equiv h_{s}^{2}.
\eeq

\begin{rem}\label{remap}
In our case, we have: $H^{k}_{s}\tilde{H}^{k}_{s}=H^{k}_{s}$.
\end{rem}

\subsection{Square root of $Q_{s}$}\label{aproot}
\beq\label{Q2}
Q(s)&=&\left(1-\frac{1}{n}\right)Id+e^{-ns}\left(\frac{1}{n}\delta_{ij}-\Theta_{i}(0)\Theta_{j}(0)\right)_{(i,j\leq n)} \nonumber \\
&\equiv& \left(1-\frac{1}{n}\right)\left(1-e^{-ns}\right)Id+e^{-ns}\cdot\sigma^{0,1}\left(\Theta(0)\right) \ .
\eeq
We are searching for $a(s)$ and $b(s)$ such that:
\beqq
\left(a(s)I+b(s)\sigma^{0,1}\left(\Theta(0)\right)\right)^{2}=Q(s),
\eeqq
or equivalently:
\beqq
\left(a(s)\right)^{2}I+2a(s)b(s)\sigma^{0,1}\left(\Theta(0)\right)+\left(b(s)\right)^{2}\sigma^{0,1}\left(\Theta(0)\right)=Q(s).
\eeqq
We compare with (\ref{Q2}) and we find:
\beqq
\left(a(s)\right)^{2}=\left(1-\frac{1}{n}\right)\left(1-e^{-ns}\right) \ \ ; \ \ 2a(s)b(s)+\left(b(s)\right)^{2}=e^{-ns}.
\eeqq
Solving this system of equations, we easily obtain:
\beq
\begin{cases}
a(s)=\sqrt{1-\frac{1}{n}}\sqrt{1-e^{-ns}} \\
b(s)=\sqrt{1-\frac{1}{n}}\left(\sqrt{1-\frac{e^{-ns}}{n-1}}-\sqrt{1-e^{-ns}}\right).
\end{cases}
\eeq

\subsection{The Ornstein-Uhlenbeck case}\label{secOU}
We consider the $n$-dimensional Ornstein-Uhlenbeck (OU) process:
\beq\label{OU}
    Z_{t} = z_{0} + \mathbb{B}_{t} - \lambda \int^{t}_{0} Z_{s} ds,
\eeq
where $\left(\mathbb{B}_{t},t\geq0\right)$ is a $n$-dimensional Brownian motion
(BM), $z_{0}\in \mathbb{R}^{n}$ and $\lambda \geq 0$.
\begin{prop}\label{propositionOU}
The Ornstein-Uhlenbeck (OU) process
$\left(\tilde{\Theta}_{t},t\geq0\right)$
in the unit sphere $S_{n-1}$ is the solution of the Stochastic Differential Equation
\beq\label{thetaOU}
\Theta^{Z}_{t}=\Theta^{Z}_{0}+\int^{t}_{0}\sigma^{0,1}(\Theta^{Z}_{s})\cdot d\hat{\mathbb{B}}_{s}-\left(\frac{n-1}{2}+\lambda\right)\int^{t}_{0}ds \ \Theta^{Z}_{s},
\eeq
where $\left(\hat{\mathbb{B}}_{t},t\geq0\right)$ is a $n$-dimensional BM.
\end{prop}

{\noindent{\textbf{Proof of Proposition \ref{propositionOU}: }}} \\
We shall study $\tilde{\varphi}_{t}\equiv\frac{Z_{t}}{|Z_{t}|}$. We remark that the Jacobi matrix
and the Hessian matrix associated respectively to the functions $\Phi(x)\equiv\frac{x}{|x|}, \ (x\neq0)$
and $g(x)\equiv |x|$ are given by:
\beqq
\left(\frac{\partial}{\partial x_{j}}\Phi_{i}(x)\right)=\frac{1}{|x|} \sigma^{0,1}(x) \ \ ; \ \
\left(\frac{\partial^{2}g(x)}{\partial x_{i}\partial x_{j}}\right)=\frac{1}{|x|} \sigma^{0,1}(x).
\eeqq
Hence, using (\ref{OU}), $\tilde{\varphi}$ satisfies the following Stochastic Differential Equation
\beq\label{phiOU}
\tilde{\varphi}_{t}&=&\tilde{\varphi}_{0}+\int^{t}_{0}\frac{1}{|Z_{s}|}\sigma^{0,1}(\tilde{\varphi}_{s})\cdot dZ_{s}-\frac{n-1}{2}\int^{t}_{0} \frac{ds}{|Z_{s}|^{2}} \tilde{\varphi}_{s} \label{phiOU2} \\
&=& \tilde{\varphi}_{0}+\int^{t}_{0}\frac{1}{|Z_{s}|}\sigma^{0,1}(\tilde{\varphi}_{s})\cdot d\mathbb{B}_{s}-\int^{t}_{0} \left(\frac{n-1}{2}\frac{\tilde{\varphi}_{s}}{|Z_{s}|^{2}}+ \frac{\lambda Z_{s}}{|Z_{s}|}\sigma^{0,1}(\tilde{\varphi}_{s}) \right) ds \nonumber \\
&=&\tilde{\varphi}_{0}+\int^{t}_{0}\frac{1}{|Z_{s}|}\sigma^{0,1}(\tilde{\varphi}_{s})\cdot d\mathbb{B}_{s}-\int^{t}_{0} \tilde{\varphi}_{s}\left(\frac{n-1}{2}\frac{1}{|Z_{s}|^{2}}+ \lambda \ \sigma^{0,1}(\tilde{\varphi}_{s}) \right) ds.
\eeq
We can replace the BM $\mathbb{B}$ by another BM $\mathbb{B}^{\ast}$:
\beqq
\mathbb{B}^{\ast}_{t}\equiv\int^{t}_{0} \left(\sigma^{0,1}(\tilde{\varphi}_{s})\cdot d\mathbb{B}_{s}+\sigma^{1,0}(\tilde{\varphi}_{s})\cdot dW_{s}\right),
\eeqq
where $\left(W_{t},t\geq0\right)$ is a BM independent from $\mathbb{B}$. \\
Thus, $\left(\gamma_{t}\equiv\int^{t}_{0}\frac{Z_{s}}{|Z_{s}|}\cdot d\mathbb{B}_{s},t\geq0\right)$ and $\mathbb{B}^{\ast}$ are two independent BMs, and from Knight's theorem (see e.g. \cite{ReY99} and other references therein), $\gamma$ is independent from the BM $\left(\hat{\mathbb{B}}_{t},t\geq0\right)$
obtained by changing the time scale of $\left(\int^{t}_{0}\frac{1}{|Z_{s}|} d\mathbb{B}^{\ast}_{s}\right)$
with the inverse of $\int^{t}_{0}\frac{ds}{|Z_{s}|^{2}}$. Finally, $\left(\Theta^{Z}_{t},t\geq0\right)$
may be obtained from $\left(\tilde{\varphi}_{t},t\geq0\right)$ by making the same change of time scale. \hfill \QED

\begin{cor}\label{corOU}
The angular part of a 2-dimensional Ornstein-Uhlenbeck process is equal to the angular part of a planar Brownian motion,
considered under the time scale $\alpha_{t}=\frac{e^{2\lambda t}-1}{2\lambda}$.
\end{cor}

\begin{rem}\label{remOU}
For further results concerning the case of a complex-valued OU process, see \cite{Vak11}.
\end{rem}

{\noindent{\textbf{Proof of Corollary \ref{corOU}: }}} \\
It follows easily from equation (\ref{OU}) for $n=2$ by taking the angular part. For the new time scale, it suffices to remark that,
with $<\cdot>$ denoting the quadratic variation of a martingale:
\beqq
\alpha_{t}\equiv <\int^{\cdot}_{0}e^{\lambda s}\cdot d\mathbb{B}_{s}>_{t}=
\int^{t}_{0} e^{2\lambda s} ds = \frac{e^{2\lambda t}-1}{2\lambda} \ .
\eeqq
\hfill \QED

\begin{rem}
Proposition \ref{propositionOU} is easily generalized for
every process of the kind:
\beq
    dZ_{s} = d\mathbb{B}_{s} + h(|Z_{s}|) Z_{s} ds,
\eeq
for every bounded function $h:\mathbb{R}^{n}\rightarrow\mathbb{R}$.
\end{rem}


\end{document}